\documentclass{article}
\usepackage{latexsym}
\usepackage{epsfig}
\usepackage{amsmath}
\usepackage{amssymb}
\usepackage{amsthm}
\newtheorem{theorem}{Theorem}[section]
\newtheorem{corollary}[theorem]{Corollary}

\newtheorem{lemma}[theorem]{Lemma}

\theoremstyle{remark}

\theoremstyle{definition}

\DeclareMathOperator\diag{diag}

\DeclareMathOperator\tr{tr}

\begin{document}

\title{Semi-definite representations for sets of cubics on the 2-sphere}

\author{Roland Hildebrand \thanks{%
Univ.\ Grenoble Alpes, CNRS, Grenoble INP, LJK, 38000 Grenoble, France
({\tt roland.hildebrand@univ-grenoble-alpes.fr}).}}

\maketitle

\begin{abstract}
The compact set of homogeneous quadratic polynomials in $n$ real variables with modulus bounded by 1 on the unit sphere $S^{n-1}$ is trivially semi-definite representable. The compact set of homogeneous ternary quartics with modulus bounded by 1 on the unit sphere $S^2$ is also semi-definite representable. This suggests that the compact set of homogeneous ternary cubics with modulus bounded by 1 on $S^2$ is semi-definite representable. We deduce an explicit semi-definite representation of this norm ball. More generally, we provide a semi-definite description of the cone of inhomogeneous ternary cubics which are nonnegative on $S^2$.
\end{abstract}

\section{Introduction}

There are few examples of spaces of polynomials in $n$ real variables of even degree such that every nonnegative polynomial is a sum of squares. These are the univariate polynomials, or equivalently, homogeneous polynomials in two variables, the quadratic polynomials, and the homogeneous ternary quartics, i.e., polynomials in three variables of degree 4 \cite{Hilbert}. In these cases the cone of nonnegative polynomials possesses a semi-definite representation \cite{NesterovSOS}, while in all other cases such a representation does not exist \cite{Scheiderer18a}. The availability of a semi-definite representation of the cone of nonnegative polynomials allows to incorporate the nonnegativity condition as a semi-definite constraint on the coefficient vector of the polynomial in a conic program.

There exist also a few examples of spaces of hermitian polynomials in complex variables where a similar result holds, namely bi-quadratic hermitian forms $p(x,y)$ with $(x,y) \in \mathbb C^2 \times \mathbb C^2$ \cite{Stormer63} and $(x,y) \in \mathbb C^2 \times \mathbb C^3$ \cite{Woronowicz}. Let ${\cal H}^n$ be the space of complex hermitian matrices of size $n$ and ${\cal H}_+^n$ the cone of positive semi-definite matrices in this space. A bi-quadratic hermitian form $p$ in $(x,y) \in \mathbb C^m \times \mathbb C^n$ is given by a matrix $H \in {\cal H}^{nm}$ such that $p(x,y) = (x \otimes y)^*H(x \otimes y)$, where $\otimes$ denotes the Kronecker product, and the asterisk denotes the complex conjugate transpose. A matrix $H \in {\cal H}^{n \cdot m}$ naturally decomposes into $m \times m$ blocks of size $n \times n$ each. Replacing each block by its complex conjugate transpose leads to the \emph{partial transpose} $H^{\Gamma}$ of the matrix. One can then establish the following sufficient semi-definite condition for nonnegativity of a form $p$:
\[ \exists\  H_1 \succeq 0,\ H_2 \succeq 0: \qquad p(x,y) \equiv (x \otimes y)^*\left(H_1 + H_2^{\Gamma}\right)(x \otimes y),
\]
where $H \succeq 0$ denotes positive semi-definiteness of the matrix $H$. Indeed, the condition implies $p(x,y) = (x \otimes y)^*H_1(x \otimes y) + (\bar x \otimes y)^*H_2(\bar x \otimes y) \geq 0$. The results in \cite{Stormer63,Woronowicz} imply that for $m = 2$, $n = 2,3$ every nonnegative bi-quadratic hermitian form can be decomposed as above, i.e., the condition is also necessary.

In fact, the statement above have been originally formulated in their dual form. In particular, we have the following result \cite{Stormer63,Woronowicz}.

{\lemma \label{lem:Woronowicz} Let $A \in {\cal H}^{2 \cdot d}$, where $d = 2$ or $d = 3$, be a complex hermitian matrix such that $A \succeq 0$, $A^{\Gamma} \succeq 0$. Then $A$ can be written as a finite sum $\sum_k (x_k \otimes y_k)(x_k \otimes y_k)^*$ of rank 1 matrices, where $x_k \in \mathbb C^2$, $y_k \in \mathbb C^d$ are non-zero vectors. \qed}

\medskip

Non-trivial cones of polynomials which are nonnegative on $\mathbb R^n$ exist only in spaces of polynomials of even degree. However, in a given space of polynomials one may also consider the closely related compact set of polynomials whose modulus does not exceed 1 on the unit sphere $S^{n-1}$, i.e., the unit ball ${\cal B}_1$ in the $\infty$-norm. This object is non-trivial independently of the parity of the degree. In a space of homogeneous polynomials of even degree $2d$, a polynomial $p(x)$ is an element of ${\cal B}_1$ if and only if the two homogeneous polynomials $\|x\|^{2d} \pm p(x)$ are nonnegative. Here $\|x\|$ is the ordinary Euclidean norm in $\mathbb R^n$.

Therefore in the space of quadratic forms in $n$ variables, the unit ball ${\cal B}_1$ is given by the semi-definite representable set $\{ H \mid -I \preceq H \preceq I \}$, where $I$ is the identity matrix. Likewise, in the space of homogeneous ternary quartics the unit ball is semi-definite representable, because the inclusion $p \in {\cal B}_1$ is equivalent to the nonnegativity of the two ternary quartics $\|x\|^4 \pm p(x)$, and the cone of nonnegative homogeneous ternary quartics is semi-definite representable.

These results suggest that the unit ball ${\cal B}_1$ in the space of homogeneous ternary cubics is also semi-definite representable, but no explicit such representation is described in the literature. In fact, a straightforward generalization of the above reasoning to spaces of polynomials of odd degree is not possible, because $\|x\|^{2d+1}$ is not a polynomial, and the use of an even power of the norm leads to the loss of homogeneity and hence the notions of nonnegativity on the unit sphere $S^{n-1}$ and on the whole space $\mathbb R^n$ become non-equivalent.

In this contribution we fill this gap by presenting a semi-definite description of the unit ball ${\cal B}_1$ in the space of homogeneous ternary cubics. This description follows from a more general result, namely a semi-definite representation of the cone of inhomogeneous ternary cubics which are nonnegative on the unit sphere $S^2$. The main idea leading to this representation is to identify the unit sphere $S^2 \subset \mathbb R^3$ with the Riemann sphere $\overline{\mathbb C}$ and to make use of Lemma \ref{lem:Woronowicz}.

Maximizing cubic polynomials over the unit sphere in $\mathbb R^n$ is an NP-hard problem \cite{NesterovMatEllipsoid} which has numerous applications in non-convex and combinatorial optimization, and several approaches for the solution of this problem have been proposed \cite{ManChoSo11,Nie12,ZhangQiYe12,AhmedStill19}. The semi-definite representations derived in this paper can be used to maximize or minimize an (in)homogeneous cubic over the 2-sphere, a task which can readily be achieved also by other methods. However, we consider a much more complicated problem than minimizing a given polynomial, namely convex optimization with constraints involving cones of nonnegative polynomials or norm balls in spaces of polynomials. Semi-definite programming can cope with such constraints only in exceptional cases, and in this contribution we discover an additional such case. An overview over descriptions (not necessarily semi-definite) of these objects in low-dimensional spaces of structured polynomials can be found in \cite{Reznick00}, see also \cite{ChoiLamReznick87} for a concrete example.

The remainder of the paper is structured as follows. In Section \ref{sec:complex} we derive a representation of polynomials on $S^2 \subset \mathbb R^3$ by complex hermitian matrices encoding hermitian rational functions on $\overline{\mathbb C}$ (Lemma \ref{lem:bijection}). In Section \ref{sec:sdr} we derive a semi-definite condition on these matrices which in the case of cubics is equivalent to nonnegativity of these rational functions (Corollary \ref{cor:dual_sdr}). In Section \ref{sec:sdr_explicit} we deduce the semi-definite representations of the cone of nonnegative inhomogeneous cubics on $S^2$ (Theorem \ref{th:nonneg_cone}) and the unit norm ball in the space of homogeneous ternary cubics (Theorem \ref{th:unit_ball}).

\section{Complex hermitian polynomials} \label{sec:complex}

Let ${\cal P}_{d,h}$ (${\cal P}_{d,ih}$) be the space of (in)homogeneous polynomials of degree $d$ in three real variables $x_1,x_2,x_3$, assembled into a vector $x \in \mathbb R^3$. The space ${\cal P}_{d,h}$ has dimension $\frac{(d+1)(d+2)}{2}$, while the dimension of ${\cal P}_{d,ih}$ equals $\sum_{j=0}^d \frac{(d+1)(d+2)}{2} = \frac{(d+1)(d+2)(d+3)}{6}$.

Every polynomial $p \in {\cal P}_{d,h}$ can be uniquely restored from knowledge of its values on the unit sphere $S^2$ by homogeneity. This is no more the case for inhomogeneous polynomials, because the polynomials in the ideal ${\cal I}$ generated by $x_1^2+x_2^2+x_3^2-1$ vanish on $S^2$. We have rather the following result.

{\lemma \label{lem:inhom_factor} Let $d > 0$ and let $p \in {\cal P}_{d,ih}$ be an arbitrary inhomogeneous polynomial. Then there exists a unique polynomial $\tilde p \in {\cal P}_{d-1,h} \oplus {\cal P}_{d,h}$ which coincides with $p$ on $S^2$. In other words, the quotient ${\cal P}_{d,ih}/{\cal I}$ is canonically isomorphic to the direct sum ${\cal P}_{d-1,h} \oplus {\cal P}_{d,h}$. }

\begin{proof}
We first construct $\tilde p$ from $p$. Let $p_{\pm}(x) = \frac{p(x) \pm p(-x)}{2}$ be the even and odd part of $p$. At every non-zero $x \in \mathbb R^3$, define the polynomial $\tilde p$ by
\[ \tilde p(x) = \left\{ \begin{array}{rcl} \|x\|^{d-1}p_+(x/\|x\|) + \|x\|^dp_-(x/\|x\|),&\quad& d\ \mbox{odd}, \\ \|x\|^{d-1}p_-(x/\|x\|) + \|x\|^dp_+(x/\|x\|),&\quad& d\ \mbox{even}, \end{array} \right.
\]
with extension to $x = 0$ by continuity. By construction we have $\tilde p \in {\cal P}_{d-1,h} \oplus {\cal P}_{d,h}$.

We now show uniqueness of $\tilde p$. Let $p \in {\cal P}_{d-1,h} \oplus {\cal P}_{d,h}$ be arbitrary. Then $p_{\pm}$ are the components of $p$ in the summands ${\cal P}_{d-1,h},{\cal P}_{d,h}$, $p_+$ corresponding to the summand consisting of even degree polynomials, and $p_-$ to the summand consisting of odd degree polynomials. If now $p(x) = 0$ for all $x \in S^2$, then both $p_{\pm}$ vanish on $S^2$ and are therefore identically zero by homogeneity. But then also $p \equiv 0$ on $\mathbb R^3$. Thus $\tilde p$ is uniquely determined by its values on $S^2$.
\end{proof}

In particular, the dimension of the quotient space ${\cal P}_{d,ih}/{\cal I}$ equals $\frac{d(d+1)}{2} + \frac{(d+1)(d+2)}{2} = (d+1)^2$. We shall refer to this quotient as the space of inhomogeneous polynomials of degree $d$ on $S^2$. 

The main result of this section is that ${\cal P}_{d,ih}/{\cal I}$ can be identified with the space ${\cal H}^{d+1}$ of complex hermitian matrices of size $d+1$. Let us index the rows and columns of matrices $H \in {\cal H}^{d+1}$ from 0 to $d$, such that, e.g., the lower right corner element of $H$ is denoted by $H_{dd}$.

\medskip

Let us first identify the unit sphere $S^2$ with the Riemann sphere $\overline{\mathbb C} = \mathbb C \cup \{\infty\}$. To every $x \in S^2$ we associate a point $z \in \overline{\mathbb C}$ by setting
\[ z(x) = \left\{ \begin{array}{rcl} \frac{-i x_2+x_3}{1 + x_1},&\quad & x_1 > -1, \\ \infty,& \quad & x_1 = -1. \end{array} \right.
\]
This entails
\begin{equation} \label{z_relations} 
\overline{z(x)} = \frac{i x_2+x_3}{1 + x_1}, \quad |z(x)|^2 = \frac{1 - x_1}{1 + x_1}, \quad 1 + |z(x)|^2 = \frac{2}{1 + x_1}
\end{equation}
for all $x \not= -e_1$, where $e_1 = (1,0,0)^T$ is the first basis vector.

For every $z \in \mathbb C \setminus \{0\}$, define the positive semi-definite rank 1 matrix
\begin{equation} \label{defZ}
Z = \frac{1}{(1 + |z|^2)^d} \begin{pmatrix} 1 \\ z \\ \vdots \\ z^d \end{pmatrix} \begin{pmatrix} 1 \\ z \\ \vdots \\ z^d \end{pmatrix}^* = \frac{1}{(1 + |z^{-1}|^2)^d} \begin{pmatrix} z^{-d} \\ \vdots \\ z^{-1} \\ 1 \end{pmatrix} \begin{pmatrix} z^{-d} \\ \vdots \\ z^{-1} \\ 1 \end{pmatrix}^* \in {\cal H}_+^{d+1},
\end{equation}
and extend this definition to $z \in \{0,\infty\}$ by continuity. The set ${\cal Z}_d = \left\{ Z \mid z \in \overline{\mathbb C} \right\}$ is a real-analytic compact manifold in the space ${\cal H}^{d+1}$. This manifold can be seen as an analog of the moment curve appearing in the study of univariate polynomials. Note that the $\mathbb R$-linear hull of ${\cal Z}_d$ is the whole space ${\cal H}^{d+1}$.

To every complex hermitian matrix $H \in {\cal H}^{d+1}$ we associate a function $p_H: S^2 \to \mathbb R$ by
\begin{equation} \label{pH} 
p_H(x) = \langle H,Z(x) \rangle,
\end{equation}
where $Z(x)$ is given by \eqref{defZ} for $z = z(x) \in \overline{\mathbb C}$, and $\langle U,V \rangle = \tr(UV)$ is the usual scalar product in the space ${\cal H}^{d+1}$.

{\lemma \label{lem:bijection} The map $H \mapsto p_H$ defines an $\mathbb R$-linear bijection between the spaces ${\cal H}^{d+1}$ and ${\cal P}_{d,ih}/{\cal I}$. }

\begin{proof}
Let us show that $p_H \in {\cal P}_{d,ih}/{\cal I}$. 
For arbitrary $0 \leq k,l \leq d$ and $x \in S^2 \setminus \{-e_1\}$ we have by virtue of \eqref{z_relations} that
\begin{align} 
\frac{z(x)^k\overline{z(x)}^l}{(1 + |z(x)|^2)^d} &= \frac{(-i x_2+x_3)^k(i x_2+x_3)^l}{2^d(1 + x_1)^{k+l-d}} = \frac{(-i x_2+x_3)^{k-\min(k,l)}(i x_2+x_3)^{l-\min(k,l)}(1 - x_1^2)^{\min(k,l)}}{2^d(1 + x_1)^{k+l-d}} \nonumber\\ &= \frac{(-i x_2+x_3)^{k-\min(k,l)}(i x_2+x_3)^{l-\min(k,l)}(1 - x_1)^{\min(k,l)}(1 + x_1)^{d-\max(k,l)}}{2^d} \label{Zx_explicit}
\end{align}
is a polynomial of degree not exceeding $d$ in $x$. Here we used that $k + l = \min(k,l) + \max(k,l)$. 

The function $p_H$ is a $\mathbb C$-linear combination of such polynomials. However, the value of $p_H$ is real by construction. Hence $p_H \in {\cal P}_{d,ih}/{\cal I}$.

$\mathbb R$-linearity of the map $H \mapsto p_H$ follows from construction. Moreover, $p_H$ is injective, because for every non-zero matrix $H \in {\cal H}^{d+1}$ the value of $p_H$ does not identically vanish on $S^2$.

But then $H \mapsto p_H$ is a bijection between ${\cal H}^{d+1}$ and ${\cal P}_{d,ih}/{\cal I}$, because the dimensions of both spaces coincide (and equal $(d+1)^2$).
\end{proof}

Denote by ${\cal C}_d$ the convex conic hull of the manifold ${\cal Z}_d$. Note that $A \succeq 0$ for every $A \in {\cal C}_d$.

By definition the matrix $H \in {\cal H}^{d+1}$ maps to a nonnegative polynomial $p_H$ on the sphere if and only if $\langle H,Z \rangle \geq 0$ for all $Z \in {\cal Z}_d$, or equivalently, $H$ is an element of the dual cone ${\cal C}^*_d$. Hence the cone of nonnegative inhomogeneous polynomials of degree $d$ on $S^2$ is isomorphic to the dual cone ${\cal C}^*_d$. Any semi-definite description of the cone ${\cal C}_d$ would then yield a semi-definite description of the cone of nonnegative polynomials on $S^2$.

\section{Semi-definite representations of ${\cal C}_d,{\cal C}^*_d$ for $d \leq 3$} \label{sec:sdr}

We shall now establish a necessary semi-definite condition on a matrix $A \in {\cal H}^{d+1}$ to be an element of ${\cal C}_d$. The main result of this section is to show that this condition is also sufficient for $d \leq 3$.

For an arbitrary matrix $A \in {\cal H}^{d+1}$, let $A_{ul},A_{ur},A_{ll},A_{lr}$ be the upper left, upper right, lower left, and lower right corner of $A$ of size $d$, respectively. Define $G_A = \begin{pmatrix} A_{ul} & A_{ur} \\ A_{ll} & A_{lr} \end{pmatrix} \in {\cal H}^{2 \cdot d}$ and its partial transpose $G_A^{\Gamma} = \begin{pmatrix} A_{ul} & A_{ll} \\ A_{ur} & A_{lr} \end{pmatrix}$. Note that the matrix $G_A$ is obtained from $A$ by duplicating the central $d-1$ rows and columns. Therefore $A \succeq 0$ is equivalent to $G_A \succeq 0$.

We have the following result.

{\lemma \label{lem:necessary_sdr} Let $A \in {\cal C}_d$. Then $A \succeq 0$, $G_A^{\Gamma} \succeq 0$. }

\begin{proof}
The relation $A \succeq 0$ follows from the definition of ${\cal C}_d$.

Let $z \in \mathbb C$ and consider the corresponding matrix $Z \in {\cal Z}_d$ defined by \eqref{defZ}. We have
\[ G_Z^{\Gamma} = \begin{pmatrix} Z_{ul} & Z_{ll} \\ Z_{ur} & Z_{lr} \end{pmatrix} = \begin{pmatrix} Z_{ul} & zZ_{ul} \\ \bar zZ_{ul} & |z|^2Z_{ul} \end{pmatrix} = \begin{pmatrix} 1 & z \\ \bar z & |z|^2 \end{pmatrix} \otimes Z_{ul} \succeq 0.
\]
For an arbitrary matrix $A \in {\cal C}_d$ the partial transpose $G_A^{\Gamma}$ lies in the convex conic hull of the manifold $\{ G_Z^{\Gamma} \mid Z \in {\cal Z}_d \} \subset {\cal H}_+^{2d}$ and is hence also positive semi-definite.
\end{proof}

{\lemma Suppose $d \leq 3$. Let $A \in {\cal H}^{d+1}$ be such that $A \succeq 0$, $G^{\Gamma}_A \succeq 0$. Then $A \in {\cal C}_d$. }

\begin{proof}
Let $d = 1$. Every positive semi-definite rank 1 matrix in ${\cal H}^2$ is in the conic hull of the manifold ${\cal Z}_1$. Any finite decomposition of $A$ into positive semi-definite rank 1 matrices then proves the claim of the lemma.

\medskip

Assume that $d = 2$ or $d = 3$. 

For every $j = 1,\dots,d-1$, let the vector $\xi_j = (0,\dots,0,1,0,\dots,0,-1,0,\dots,0)^T \in \mathbb C^{2d}$ be such that its $(j+1)$-th element equals 1, its $(d+j)$-th element equals $-1$, and all other elements equal zero. Then $G_A\xi_j = 0$ by construction of $G_A$.

By the assumptions of the lemma the matrix $G_A$ satisfies the conditions of Lemma \ref{lem:Woronowicz}. Hence $G_A$ can be decomposed into a finite sum $\sum_k (w_k \otimes y_k)(w_k \otimes y_k)^*$ of positive semi-definite rank 1 matrices, where $w_k \in \mathbb C^2$, $y_k \in \mathbb C^d$ are non-zero vectors. The vectors $\xi_j$ are in the kernel of every of these rank 1 matrices, which entails $\xi_j^*(w_k \otimes y_k) = (w_k)_1(y_k)_{j+1} - (w_k)_2(y_k)_j = 0$ for all $k$ and all $j = 1,\dots,d-1$. For every $k$, form the vector 
\[ \psi_k = \left((w_k)_1(y_k)_1,(w_k)_2(y_k)_1 = (w_k)_1(y_k)_2,\dots,(w_k)_2(y_k)_{d-1} = (w_k)_1(y_k)_d,(w_k)_2(y_k)_d\right)^T \in \mathbb C^{d+1},
\]
constructed such that the Kronecker product $w_k \otimes y_k$ is obtained from $\psi_k$ by duplicating the central $d-1$ elements. Therefore $G_A = \sum_k (w_k \otimes y_k)(w_k \otimes y_k)^*$ implies $A = \sum_k \psi_k\psi_k^*$. Moreover, the rank 1 matrix $\psi_k\psi_k^*$ is a multiple of the matrix $Z \in {\cal Z}_d$ corresponding to $z = \frac{(w_k)_2}{(w_k)_1} \in \overline{\mathbb C}$. But this entails $A \in {\cal C}_d$.
\end{proof}

{\theorem \label{th:sdr_equivalence} Suppose $d \leq 3$. Then a matrix $A \in {\cal H}^{d+1}$ is an element of the cone ${\cal C}_d$ if and only if $A \succeq 0$, $G_A^{\Gamma} \succeq 0$. }

\begin{proof}
The theorem follows by applying the two preceding lemmas.
\end{proof}

The conditions $A \succeq 0$, $G^{\Gamma}_A \succeq 0$ are $\mathbb R$-linear in $A$ and define semi-definite constraints on $A$. This yields a semi-definite description of the cone ${\cal C}_d$ for $d \leq 3$. In particular, ${\cal C}_d$ is a spectrahedral cone, namely the preimage of the positive semi-definite matrix cone ${\cal H}_+^{3d+1}$ under the injective $\mathbb R$-linear map $L: {\cal H}^{d+1} \to {\cal H}^{3d+1}$ defined by $L: A \mapsto \diag\left(A,G_A^{\Gamma}\right)$.

{\corollary \label{cor:dual_sdr} Suppose $d \leq 3$. Then a matrix $H \in {\cal H}^{d+1}$ is an element of the dual cone ${\cal C}^*_d$ if and only if there exist $B \in {\cal H}_+^{d+1}$, $C \in {\cal H}_+^{2d}$ such that for every $A \in {\cal H}^{d+1}$ we have
\[ \left\langle \diag(B,C),\diag\left(A,G_A^{\Gamma}\right) \right\rangle = \langle H,A \rangle.
\] }

\begin{proof}
The corollary is obtained from Theorem \ref{th:sdr_equivalence} by conic duality.
\end{proof}

For general degree $d$ the semi-definite conditions in Corollary \ref{cor:dual_sdr} are only sufficient for the inclusion $H \in {\cal C}^*_d$, by applying conic duality to the assertion of Lemma \ref{lem:necessary_sdr}.

\section{Explicit semi-definite representations} \label{sec:sdr_explicit}

In this section we shall develop the semi-definite constraints in Corollary \ref{cor:dual_sdr} explicitly for $d = 3$. This allows to obtain a semi-definite representation of the unit norm ball ${\cal B}_1$ in the space of homogeneous ternary cubics.

With $z = z(x)$ the matrix in \eqref{defZ} is by virtue of \eqref{Zx_explicit} given by
\[ Z(x) = \frac18 \begin{pmatrix} (1 + x_1)^3 & (i x_2+x_3)(1 + x_1)^2 & (i x_2+x_3)^2(1 + x_1) & (i x_2+x_3)^3 \\ (-i x_2+x_3)(1 + x_1)^2 & (1 - x_1)(1 + x_1)^2 & (i x_2+x_3)(1 - x_1^2) & (i x_2+x_3)^2(1 - x_1) \\ (-i x_2+x_3)^2(1 + x_1) & (-i x_2+x_3)(1 - x_1^2) & (1 - x_1)^2(1 + x_1) & (i x_2+x_3)(1 - x_1)^2 \\ (-i x_2+x_3)^3 & (-i x_2+x_3)^2(1 - x_1) & (-i x_2+x_3)(1 - x_1)^2 & (1 - x_1)^3 \end{pmatrix}.
\]

We have the following characterization of nonnegative cubic polynomials.

{\theorem \label{th:nonneg_cone} The polynomial $p = \sum_{2 \leq j+k+l \leq 3} c_{jkl}x_1^jx_2^kx_3^l \in {\cal P}_{2,h} \oplus {\cal P}_{3,h}$ is nonnegative on $S^2$ if and only if there exist $B \in {\cal H}_+^4$, $C \in {\cal H}_+^6$ such that
\begin{align*}
B_{11} + C_{11} &= c_{200} + c_{300} \\
B_{44} + C_{66} &= c_{200} - c_{300} \\
B_{12} + C_{12} + C_{41} &= c_{201} + c_{101} + i(c_{210} + c_{110}) \\
B_{34} + C_{63} + C_{56} &= c_{201} - c_{101} + i(c_{210} - c_{110}) \\
B_{22} + C_{22} + C_{44} + C_{15} + C_{51} &= 2c_{020} + 2c_{002} - c_{200} + 2c_{120} + 2c_{102} - 3c_{300} \\
B_{33} + C_{33} + C_{55} + C_{26} + C_{62} &= 2c_{020} + 2c_{002} - c_{200} - 2c_{120} - 2c_{102} + 3c_{300} \\
B_{23} + C_{16} + C_{23} + C_{52} + C_{45} &= 3c_{003} + c_{021} - 2c_{201} + i(3c_{030} + c_{012} - 2c_{210}) \\
B_{41} + C_{34} &= c_{003} - c_{021} + i(c_{030} - c_{012}) \\
B_{24} + C_{53} + C_{46} &= c_{002} - c_{020} + c_{120} - c_{102} + i(c_{011} - c_{111}) \\
B_{13} + C_{13} + C_{42} &= c_{002} - c_{020} - c_{120} + c_{102} + i(c_{011} + c_{111}).
\end{align*} }

\begin{proof}
Let $H \in {\cal H}^4$ be the matrix corresponding to the polynomial $p = p_H$ as in \eqref{pH}. Then $p$ is nonnegative on $S^2$ if and only if $H \in {\cal C}^*_d$, i.e., there exist matrices $B \in {\cal H}_+^4$, $C \in {\cal H}_+^6$ satisfying the linear conditions in Corollary \ref{cor:dual_sdr}. Since the matrices $Z(x)$ span the whole space ${\cal H}^4$ when $x$ runs through $S^2$, these conditions are equivalent to
\[ \left\langle \diag(B,C),\diag\left(Z(x),G_{Z(x)}^{\Gamma}\right) \right\rangle = \langle H,Z(x) \rangle = p(x), \qquad \forall\ x \in S^2.
\]
Comparing the coefficients of the polynomials on the left- and right-hand side modulo $\|x\|^2 - 1$ yields the equations in the theorem.
\end{proof}

Let now $p \in {\cal P}_{3,h}$ be a homogeneous ternary cubic. Then the modulus of $p$ is bounded by 1 on the sphere $S^2$ if and only if the inhomogeneous cubic $\|x\|^2 + p(x)$ is nonnegative on $S^2$. This yields the following semi-definite representation of ${\cal B}_1$.

{\theorem \label{th:unit_ball} The homogeneous ternary cubic $p = \sum_{j+k+l = 3} c_{jkl}x_1^jx_2^kx_3^l \in {\cal P}_{3,h}$ has a modulus bounded by 1 on $S^2$ if and only if there exist $B \in {\cal H}_+^4$, $C \in {\cal H}_+^6$ such that
\begin{align*}
B_{11} + C_{11} &= 1 + c_{300} \\
B_{44} + C_{66} &= 1 - c_{300} \\
B_{12} + C_{12} + C_{41} = B_{34} + C_{63} + C_{56} &= c_{201} + ic_{210} \\
B_{22} + C_{22} + C_{44} + C_{15} + C_{51} &= 3 + 2c_{120} + 2c_{102} - 3c_{300} \\
B_{33} + C_{33} + C_{55} + C_{26} + C_{62} &= 3 - 2c_{120} - 2c_{102} + 3c_{300} \\
B_{23} + C_{16} + C_{23} + C_{52} + C_{45} &= 3c_{003} + c_{021} - 2c_{201} + i(3c_{030} + c_{012} - 2c_{210}) \\
B_{41} + C_{34} &= c_{003} - c_{021} + i(c_{030} - c_{012}) \\
B_{24} + C_{53} + C_{46} = -(B_{13} + C_{13} + C_{42}) &= c_{120} - c_{102} - ic_{111}.
\end{align*} }

\begin{proof}
The result follows from Theorem \ref{th:nonneg_cone} by setting $c_{200} = c_{020} = c_{002} = 1$, $c_{110} = c_{101} = c_{011} = 0$.
\end{proof}

\bibliographystyle{plain}
\bibliography{convexity,polynomial_optim}

\end{document}